\documentclass[12pt]{amsart}
\usepackage{amsmath}
\usepackage{amsfonts}
\theoremstyle{plain}
\newtheorem{thm}{Theorem}[section]
\newtheorem{lmm}[thm]{Lemma}

\theoremstyle{definition}

\newtheorem{rmk}[thm]{Remark}

\def\op{\operatorname}

\def\un{\underline}

\begin{document}

\title{Rigidity of the Minimal Grope Group}
\author{Matija Cencelj}
\address{IMFM, University of Ljubljana, 
Jadranska 19, Ljubljana, SLOVENIA}
\email{matija.cencelj@guest.arnes.si}
\author{Katsuya Eda}
\address{School of Science and Engineering, 
Waseda University, Tokyo 169-8555, JAPAN}
\email{eda@waseda.jp}
\author{Ale\v s Vavpeti\v c}
\address{IMFM, University of Ljubljana, 
Jadranska 19, Ljubljana, SLOVENIA}
\email{ales.vavpetic@fmf.uni-lj.si}

\date{August 3, 2006}


\thanks{Supported in part by the Slovenian-Japanese research grant BI--JP/05-06/2,
ARRS research program No. 0101-509, the ARRS
research project of Slovenia No. J1--6128--0101--04 and 
the Grant-in-Aid for Scientific research (C) of Japan 
No. 16540125.} 

\maketitle

\begin{abstract}
{We give a systematic definition of the fundamental groups of gropes, which we call
grope groups. We show that there exists a nontrivial homomorphism from the minimal grope group $M$
to another grope group $G$ only if $G$ is the free product
of $M$ with another grope group.}
\end{abstract}

\section{Introduction}
Here we study groups whose classifying spaces are (open infinite) gropes (a recent
short note on gropes in general is \cite{Teichner}).
In algebra these groups first
appeared in the proof of a lemma by Alex Heller \cite{Heller} as follows. Let
$\varphi_0$ be a homomorphism  from the free group $F_0$ on one generator $\alpha$
to any perfect group $P$. Let 
$$
\varphi_0(\alpha)=[p_0,p_1][p_2,p_3]\cdots[p_{2n-2},p_{2n-1}]\in P 
\qquad (*) 
$$ 
then
we can extend $\varphi_0$ to a homomorphism $\varphi_1$ of a
(nonabelian) free group $F_1$ on $2n$ generators
$\beta_0,\ldots,\beta_{2n-1}$ by setting $\varphi_1(\beta_i)=p_i$. Note
that $\varphi_0(\alpha_0)$ may have several different expressions as a
product of commutators, so we may choose any; even if some of the
elements $p_1,\ldots,p_{2n-1}$ coincide we
let all elements $\beta_i$ to be distinct. Now we repeat the above
construction for every homomorphism $\varphi_1|_{\langle
\beta_i\rangle}$ of the free group on one generator to $P$ and thus
obtain a homomorphism $\varphi_{2}:F_{2}\rightarrow P$.
Repeating the above construction we obtain a direct system of inclusions
of free groups $F_1\rightarrow F_2\rightarrow F_3\rightarrow 
\cdots$ and homomorphisms $\varphi _n:F_n\rightarrow P$. The direct 
limit of $F_n$ is a locally free perfect group $D$ and every group 
obtained by the above construction is called a grope group (and its
clasifying space is a grope). This construction shows therefore 
that every homomorphism from a free group on one generator to a perfect group $P$
can be extended to a homomorphism from a grope group to $P$. Note that in case the
perfect group $P$ has the Ore property (\cite{Ore}, \cite{EllersGordeev}) that every
element in $P$ is a commutator, in the above process ($*$)
we can make every generator in the chosen basis of $F_n$ a single commutator of two
basis elements of $F_{n+1}$. The group
obtained in this way is the minimal grope group $M$. Clearly every grope group admits
many epimorphisms onto $M$. In the sequel we show that $M$ admits a nontrivial homomorphism to another grope group $G$
only if the latter is the free product $G\cong M*K$ where $K$ is a grope group. 

Gropes were introduced by \v Stan'ko \cite{Stanko}. They have an important role in
geometric topology 
(\cite{Cannon}, for more recent use in dimension theory see
\cite{DranishnikovRepovs} and \cite{CenceljRepovs}). Their fundamental
groups were used by Berrick and Casacuberta to show that the plus-construction in
algebraic K-theory is localization 
\cite{BerrickCasacuberta}. Recently \cite{BadziochFeshbach} such
a group has appeared in the construction of a perfect group with a nonperfect
localization.

In the first part of the paper we give a systematic definition of grope groups and
prove some technical lemmas. In the second part
we prove that the minimal grope group admits nontrivial homomorphisms to almost no
other grope group thus proving that there exist at
least two distinct grope groups.

\section{Systematic definition of grope groups and basic facts}
For every positive integer $n$ let $\underline{n}=\{0,1,\ldots,n-1\}$.
The set of non-negative integers is denoted by $\mathbb{N}$. 
We denote the set of finite sequences of elements of a set $X$ by
$Seq(X)$ and the length of a sequence $s\in Seq(X)$ by $lh(s)$. The
empty sequence is denoted by $\emptyset$. 

For a non-empty set $A$ let $L(A)$ be the set $\{ a,a^-: a\in A\}$,
which we call the set of letters. We identify $(a^-)^-$ with $a$. 
Let $\mathcal{W}(A) = Seq(L(A))$, which we call the set of words. 
For a word $W \equiv a_0\cdots a_n$, define $W^- \equiv a_n^-\cdots
a_0^-$. We write $W\equiv W'$ for identity in $\mathcal{W}(A)$
while $W=W'$ for identity in the free group generated by $A$. For
instance $aa^- = \emptyset$ but $aa^- \not\equiv \emptyset$. 
 We adopt $[a,b] = aba^{-1}b^{-1}$ as the definition of a commutator.

To describe all the grope groups we introduce some notation. 

A {\it grope frame} $S$ is a subset of $Seq(\mathbb{N})$ satisfying: 
$\emptyset \in S$ and for every $s\in S$ there exists $n>0$ such that
$\underline{2n} = \{ i\in \mathbb{N}: si \in S\}$. 

For each grope frame $S$ we induce formal symbols $c^S_s$ for $s\in S$
and define $E^S_m =\{ c^S_s: lh(s) = m, s\in S\}$ and a free group
$F^S_m = \langle E^S_m\rangle$. 
Then define $e^S_m:F^S_m \to F^S_{m+1}$ by: 
$e^S_m(c^S_s) = [c_{s0}^S,c_{s1}^S]\cdots [c_{s\, 2k-2}^S,c_{s\, 2k-1}^S]$
where $\un{2k} =\{ i\in \mathbb{N} : si \in S\}$. 
Let $G^S = \varinjlim(F^S_m,e^S_m: m\in \mathbb{N} )$ and $e^S_{mn} =
e^S_{n-1}\cdots e^S_m$ for $m\le n$ and every such group $G^S$ is a grope group. 

For $s\in S$, $s$ is {\it binary branched}, if $\{ i\in \mathbb{N} : si\in S\}
= \un{2}$. 
Let $S_0$ be a grope frame such that every $s\in S_0$ is binary
branched, i.e. $S_0 = Seq(\un{2})$. Then $G^{S_0}=M$ is the so-called minimal
grope group. 
Since $e^S_m$ is injective, we frequently regard $F^S_m$ is a subgroup
of $G^S$. 

For a non-empty word $W$ the {\it head} of $W$ is the left most letter $b$ of
$W$, i.e. $W\equiv bX$ for some word $X$, and  the {\it tail} of $W$ is
the right most letter $c$ of $W$, i.e. $W\equiv Yc$ for some word $Y$. 
When $AB\equiv W$, we say that $A$ is the head part of $W$ and $B$ is
the tail part of $W$. 
For a word $W\in \mathcal{W}(E_m^S)$ and $n\ge m$, we let $e_{m, n}^S[W]$ be
a word in $\mathcal{W}(E_n^S)$ defined as follows: 
$ e_{m,m}^S[W] \equiv W$ and 
$e_{m, n+1}^S[W]$ is obtained by replacing every $c_t$ in $e_{m,n}^S[W]$ by 
\[
 \op{(P0)}\quad c_{t0}^Sc_{t1}^Sc_{t0}^{S-}c_{t1}^{S-}\cdots c_{t\, 2k-2}^Sc_{t\,
 2k-1}^Sc_{t\, 2k-2}^{S-}c_{t\, 2k-1}^{S-}
\]
and every $c_t^{S-}$ by 
\[
 \op{(P1)}\quad c_{t\, 2k-1}^S c_{t\, 2k-2}^Sc_{t\, 2k-1}^{S-}c_{t\, 2k-2}^{S-}
 \cdots c_{t1}^Sc_{t0}^Sc_{t1}^{S-}c_{t0}^{S-}
\]
respectively. 

We drop the superscript ${}^S$, if no confusion can occur. 

For a reduced word $W\in \mathcal{W}(E_n)$ with $W\in F_m$ for $m<n$, 
let $W_0\in \mathcal{W}(E_m)$ such that $e_{m,n}[W_0] \equiv W$. (The
existence of $W_0$ is assured in Lemma~\ref{lmm:basic2}.) 
A subword $V$ of $W$ is {\it small}, if there exists a letter $c_s$ or
$c_s^-$ in $W_0$ and $i\in \mathbb{N}$ such that $V$ is a subword of $e_{m+1,
n}[c_{si}]$
or $e_{m+1, n}[c_{si}^-]$ respectively.
(Note that being small depends on $m$. In the following usage of this
notion $m$ and $n$ are always fixed in advance.)

\medskip

\noindent
{\bf Observation 1.}
Let $n>m+1$ and let $W\equiv e_{m+1, n}[c_{s0}]$. Suppose that $X\in
\mathcal{W}(E_n)$ is a reduced word and $X\in F_m$. When $W$ is a
subword of $X$, $W$ may appear in $e_{m,n}[c_s]$ or $e_{m,n}[c_s^-]$ and
hence we cannot uniquely determine a successive letter to $W$ in
$X$. However, if $X\equiv WY$ for some $Y$, the head of $Y$ is uniquely
determined as $c_{s10\cdots 0}$. Also if we know the preceding letter to
$W$, i.e. $X\equiv ZWY$ and we know the tail of $Z$ which is
$c_{t10\cdots 0}^-$ or $c_{t0\cdots 0}^-$ for some $t$ or $c_{s110\cdots
0}^-$, the head of $Y$ is uniquely determined. 
That is, the head is $c_{s10\cdots 0}$, if the preceding letter is
$c_{t10\cdots 0}^-$ or $c_{t0\cdots 0}^-$ for some $t$ and the head is 
$c_{s110\cdots 0}$, if the preceding letter is $c_{s110\cdots 0}^-$. 
(In the above we ignore $i$-digits for $i>n$.)

\medskip

\noindent
{\bf Observation 2.}
A letter $c_{s0\cdots 0}\in
\mathcal{W}(E_n)$ for $lh(s) = m$ possibly appears in $e_{m,n}[W_0]$ in
the following cases. 
When $n = m+1$, 
$c_{s0}$ appears once in $e_{m,n}[c_s]$ and also once in $e_{m,n}[c_s^-]$. 
According to the increase of $n$, $c_{s0\cdots 0}$ appears in many
parts. 
$c_{s0\cdots 0}$ appears $2^{n-m-1}$-times in $e_{m,n}[c_s]$ and also 
$2^{n-m-1}$-times in $e_{m,n}[c_s^-]$. 
\begin{lmm}\label{lmm:basic1}
For a word $W\in \mathcal{W}(E_m)$ and $n\ge m$, 
$e_{m,n}[W]$ is reduced, if and only if $W$ is reduced.
\end{lmm}
\begin{lmm}\label{lmm:basic2}
For a reduced word $V\in \mathcal{W}(E_n)$ and $n\ge m$, 
$V\in F_m$ if and only if there exists $W\in \mathcal{W}(E_m)$ such
 that $e_{m,n}[W] \equiv V$.
\end{lmm}
\begin{proof}
The sufficiency is obvious. To see the other direction, let $W$ be a
 reduced word in $\mathcal{W}(E_m)$ such that $e_{m,n}[W] = V$ in $F_n$. By
 Lemma~\ref{lmm:basic1} $e_{m,n}[W]$ is reduced. Since every element in
 $F_n$ has a unique reduced word in $\mathcal{W}(E_n)$ presenting
 itself, we have $e_{m,n}[W] \equiv V$.
\end{proof}
\begin{lmm}\label{lmm:basic3}
Let $m<n$ and $A$ be a non-empty word in $\mathcal{W}(E_n)$. 
Let $X_0AY_0$ and $X_1AY_1$ be reduced words in $\mathcal{W}(E_n)$
 satisfying $X_0AY_0, X_1AY_1\in F_m$. 
\begin{itemize}
\item[(1)] If $A$ is not small, $X_0A\notin F_m$ and $X_1A\notin F_m$,
           then the heads of $Y_0$ and $Y_1$ are the same. 
\item[(2)] Let $X_0$ be an empty word. If $A$ is not small and 
           $A\notin F_m$, the heads of $Y_0$ and $Y_1$ are the same.
\item[(3)] Let $X_0$ and $X_1$ be empty words. If $A\notin F_m$, the
           heads of $Y_0$ and $Y_1$ are the same.
\end{itemize} 
\end{lmm}
\begin{proof}
(1) Since $X_0AY_0 \in F_m$ but $X_0A\notin F_m$, we have a letter
 $c\in E_m\cup E_m^-$ and words $U_0,U_1,U_2$ such that 
$U_1\not\equiv \emptyset$, $U_2\not\equiv \emptyset$, $X_0A\equiv U_0U_1$
 and $U_1U_2\equiv e_{m,n}[c]$.  
Since $A$ is not small, $c$ and $U_0,U_1,U_2$ are uniquely determined by
 $A$. Since the same thing holds for $X_1AY_1$, we have the conclusion
 by Observation 1 for $n>m+1$. (The case for $n=m+1$ is easier.) 

(2) Since $AY_0\in F_m$, $A\notin F_m$ and $A$ is not a small word, for any word $B$
 such that $BA$ is reduced we have $BA\notin F_m$. In particular $X_1A\notin F_m$
 and the conclusion follows from (1). 

(3) Since $AY_0\in F_m$, there are $A_0$ and non-empty
 $U_0,U_1$ such that $A_0\in F_m$, $A\equiv A_0U_0$ and $U_0U_1 \equiv 
e_{m,n}[c]$ for some $c\in E_m\cup E_m^-$. 
Since $A\notin F_m$, the head of $U_1$ is uniquely determined by $A$ and
 hence the heads of $Y_0$ and $Y_1$ are the same (Observation 1). 
\end{proof}
\begin{lmm}\label{lmm:basic5}
Let $m<n$ and $A,X,Y$ in $\mathcal{W}(E_n)$ and $AXA^-Y\in F_m$. 
If $AXA^-Y$ is reduced and $A$ is not small, then $AXA^-\in F_m$ and
 $Y\in F_m$. 
\end{lmm}
\begin{proof}
The head of the reduced word in $\mathcal{W}(E_m)$ for the element
 $AXA^-Y$ is $c_s$ or $c_s^-$ for $c_s\in E_m$. According to $c_s$ or
 $c_s^-$, $A\equiv e_{m+1, n}[c_{s0}]Z$ or $e_{m+1, n}[c_{sk}]Z$ 
for a non-empty word $Z$, where $\un{k+1} = \{ i\in \mathbb{N} : si\in S\}$
 is even. 
Then $A^- \equiv Z^-e_{m+1, n}[c_{s0}^-]$ or 
$A^- \equiv Z^-e_{m+1, n}[c_{sk}^-]$ and hence $AXA^-\in F_m$ and
 consequently $Y\in F_m$. 
\end{proof}
\begin{lmm}\label{lmm:antinormal}
For $e\neq x\in F_m^S$ and $u\in G^S$, 
$uxu^{-1}\in F_m^S$ implies $u\in F_m^S$. 
\end{lmm}
\begin{proof}
There exists $n\ge m$ such that $u\in F_n$.
Let $W$ be a cyclically reduced word and $V$ be a reduced word such
 that $x = VWV^-$ in $F_m$ and $VWV^-$ is reduced. Then 
$e_{m,n}(x) = e_{m,n}[V]e_{m,n}[W]e_{m,n}[V]^-$ and 
$e_{m,n}[V]$ is reduced and $e_{m,n}[W]$ is cyclically reduced by
 Lemma~\ref{lmm:basic1}. Let $U$ be a reduced word for $u$ in $F_n$. 
Let $k = lh(U)$. 
Then $e_{m,n}(x^{2k+1}) = e_{m,n}[V]e_{m,n}[W]^{2k+1}e_{m,n}[V]^-$ and the
 right hand term is a reduced word. 
Hence the reduced word for $ux^ku^-$ of the form $Xe_{m,n}[W]Y$, where 
$Ue_{m,n}[V]e_{m,n}[W]^k = X$ and $e_{m,n}[W]^ke_{m,n}[V]^-U^- = Y$. Since
 $ux^ku^{-1} \in F_m$, $X\in F_m$ and $Y\in F_m$. Now we have
 $Ue_{m,n}[V]\in e_{m,n}(F_m)$ and hence $U\in e_{m,n}(F_m)$, which
 implies the conclusion. 
\end{proof}
\begin{lmm}\label{lmm:reducedform}
Let $UWU^-$ be a reduced word in $\mathcal{W}(E_n)$. 
If $UWU^-\in F_m$ and $W$ is cyclically reduced, then $U, W \in F_m$. 
\end{lmm}
\begin{proof}
If $U$ is empty or $n=m$, then the conclusion is obvious. If $U\in F_m$, then
 $WU^-\in F_m$ and so $W\in F_m$. Suppose that $U$ is $U\not\in F_m$. Since 
$UWU^-, UW^-U^-\in F_m$, the head of $W$ and that of $W^-$ is the same
 by Lemma~\ref{lmm:basic3} (3), which contradicts that $W$ is cyclically
 reduced. 
\end{proof}
\begin{lmm}\label{lmm:transform}
Let $XY$ and $YX$ be reduced words in $\mathcal{W}(E_n)$ for 
$n\ge m$. If $XY$ and $YX$ belong to $F_m$, then both of $X$ and
 $Y$ belong to $F_m$.   
\end{lmm}
\begin{proof}
We may assume $n>m$. 
When $n>m$, the head of $e_{m,n}[W]$ for a non-empty word 
$W\in \mathcal{W}(E_m)$ is $c_{s0\cdots 0}$ or $c_{sk0\cdots 0}$ where
 $lh(s)= m$ and $k+1 = \{ i\in \mathbb{N} : si\in S\}$ is even. (When $n=m+1$, there
 appears no $0\cdots 0$.) 
Since $X^-Y^-\in F_m$ and $X^-Y^-$ is reduced, the tail of $X$ is of the
 form $c_{s0\cdots 0}^-$ or $c_{sk0\cdots 0}^-$. We only deal with the
 former case. 
Suppose that $X\notin F_m$. Since $XY\in F_m$ and $XY$ is reduced, 
 $X \equiv Ze_{m+1, n}[c_{s1}c_{s0}^-]$ for some $Z$. This implies 
$X^-\equiv e_{m+1, n}[c_{s0}c_{s1}^-]Z^-$, which contradicts that
 $X^-Y^-\in F_m$ and $X^-Y^-$ is reduced. 
Now we have $X,Y\in F_m$. 
\end{proof}
\begin{lmm}\label{lmm:three}
Let $m<n$ and $A,B,C$ in $\mathcal{W}(E_n)$ and $e\neq ABCA^-B^-C^-\in F_m$. 
If $ABCA^-B^-C^-$ is a reduced word and at least one of $A,B,C$ is not
 small, then $A,B,C\in F_m$. 
\end{lmm}
\begin{proof}
Since $ABCA^-B^-C^-\neq e$, at most one of $A,B,C$ is empty. 
When $C$ is empty, the conclusion follows from Lemma~\ref{lmm:basic5}
 and the fact that $BAB^-A^-$ is also reduced and $BAB^-A^-\in F_m$. 

Now we assume that $A,B,C$ are non-empty. If $A$ is not small, then
 $ABCA^-\in F_m$ and $B^-C^-\in F_m$ by Lemma~\ref{lmm:basic5}. Since
 $BC$ is cyclically reduced, $A\in F_m$ and $BC\in F_m$ by
 Lemma~\ref{lmm:reducedform}. The conclusion follows from
 Lemma~\ref{lmm:transform}. In the case that $C$ is not small, the
 argument is similar. The remaining case is when $A$ and $C$ are small. 
Then $ABCA^-B^-C^-\in F_m$ and $CBAC^-B^-A^-\in F_m$ imply $A\equiv C$,
 which contradicts the reducedness of $ABCA^-B^-C^-$. 
 \end{proof}

\begin{lmm}\label{lmm:empty}
Let $m<n$ and $A,B,C$ in $\mathcal{W}(E_n)$ and $e\neq 
ABCA^-B^-C^-\in F_m$. If $ABCA^-B^-C^-$ is a reduced word and $A,B,C$
 are small, then one of $A,B,C$ is empty.
 
 Assume $C$ is empty. Then there exists $c_s\in E_m$ such that $s$ is binary
branched and 
either 

\begin{center}
$A\equiv e_{m+1, n}[c_{s0}]$ and 
$B\equiv e_{m+1, n}[c_{s1}]$,
\end{center}
 or
\begin{center}
$A\equiv e_{m+1, n}[c_{s1}]$ and $B\equiv e_{m+1, n}[c_{s0}]$. 
\end{center}
\end{lmm}

\begin{proof}
Since $A,B,C$ are small, all the words $A,B,C$ and their inverses must 
 be subwords of $e_{m+1, n}[c_{si}]$, $i=0,1$, or $e_{m+1, n}[c_{si}^-]$, 
for an element $c_s\in E_m$, and
in particular that either
\[
 ABCA^-B^-C^- = e_{m,n}(c_s) = e_{m+1, n}[c_{s0}c_{s1}c_{s0}^-c_{s1}^-]
\] 
or 
\[
 ABCA^-B^-C^- = e_{m,n}(c_s^-) = e_{m+1, n}[c_{s1}c_{s0}c_{s1}^-c_{s0}^-],
\]
where the left most and right most terms are reduced words. 
We remark that if the cardinality of 
$\{ i\in \mathbb{N} : si\in S\}$ were greater than $2$, one of $A,B,C$ would not be
small; hence
in our case $s$ is
 binary branched. 

We only deal with the first case. Then 
$ABC\equiv e_{m+1, n}[c_{s0}c_{s1}]$ and 
$A^-B^-C^- \equiv e_{m+1, n}[c_{s0}^-c_{s1}^-]$. In case $A,B,C$ are
 non-empty, $A$ is a proper subword of $e_{m+1, n}[c_{s0}]$ or 
$C$ is a proper subword of $e_{m+1, n}[c_{s1}]$. In either case  
$A^-B^-C^- \equiv e_{m+1, n}[c_{s0}^-c_{s1}^-]$ does not hold. Hence
 one of $A,B,C$ is empty. We may assume $C$ is empty. Since $A,B$ are
 small, $A\equiv e_{m+1, n}[c_{s0}]$ and $B\equiv e_{m+1, n}[c_{s1}]$. 
\end{proof}%

\section{Proof of Theorem~\ref{thm:main}}
In this section we prove 
\begin{thm}\label{thm:main}
The minimal grope group $M=G^{S_0}$ admits a nontrivial homomorphism into a grope group
 $G^S$, if and only if there exists $s\in S$ such that 
a frame $\{ t\in Seq(\mathbb{N}): st\in S\}$ is equal to $S_0$. 
\end{thm}

It is easy to see that the condition on $G^S$ in the above theorem is
equivalent to $G^S\cong M*K$, where $K$ is another grope group.

In our proof of Lemma~\ref{lmm:main} we analyze a reduction procedure
of a word $Y^-ABYX^-A^-B^-X$ where $Y^-ABY$ and $X^-A^-B^-X$ are
reduced. Lemmas~\ref{lmm:induction1}, \ref{lmm:induction2},
\ref{lmm:induction3} and \ref{lmm:induction4} show connections between
our reduction steps in case at least one of $X$ and $Y$ is empty. 
Lemma~\ref{lmm:crucial1} corresponds to the final step, i.e. when we
have the reduced word. Lemmas~\ref{lmm:crucial2} and
\ref{lmm:remaining} correspond to the case that $X$ and $Y$ are
non-empty. In the following lemmas we assume $m < n$. 

\begin{lmm}\label{lmm:induction1}
Let $A,B\in \mathcal{W}(E_n)$ be non-empty reduced words such
 that $ABA^-B^- \neq e$ and $AB$ and $A^-B^-$ are reduced words.
Then the following hold: 
\begin{itemize}
\item[(1.1)] If $B\equiv B_0A$, then $B_0$ is non-empty, $AB_0$ and
             $A^-B_0^- $ are reduced words and $AB_0A^-B_0^- = ABA^-B^-$. 
             In addition if $AB_0,A^-B_0^-\in F_m$, then 
             $AB,A^-B^-\in F_m$. 
\item[(1.2)] If $A\equiv A_0B$, then $A_0$ is non-empty, $A_0B$ and
             $A_0^-B^-$ are reduced words and $A_0BA_0^-B^- 
             = ABA^-B^-$. In addition if $A_0B,A_0^-B^-\in F_m$, 
             then $AB,A^-B^-\in F_m$. 
\item[(1.3)] If $A\equiv A_0Z$ and $B\equiv B_0Z$ for non-empty words
             $A_0$ and $B_0$ and $B_0A_0^-$ is reduced, then
             $A_0ZB_0A_0^-Z^-B_0^-$ is reduced and 
             $A_0ZB_0A_0^-Z^-B_0^- = ABA^-B^-$. In addition if 
             $A_0, B_0, Z\in F_m$, then $AB,A^-B^-\in F_m$. 
\end{itemize}
\end{lmm}
\begin{proof}
We only show (1.1). The non-emptiness of $B_0$ follows from
 $ABA^-B^-\neq e$. Since $AB$ and $A^-B^-$ are reduced, $AB_0$ and
 $A^-B_0^-$ are cyclically reduced and hence the second statement
 follows from Lemma~\ref{lmm:transform}. 
\end{proof}
\begin{lmm}\label{lmm:induction2}
Let $A,B,C\in \mathcal{W}(E_n)$ be reduced words (possibly empty) such
 that $ABCA^-B^-C^- \neq e$ and $AB$ and $CA^-B^-C^- $ are reduced words.
Then the following hold: 
\begin{itemize}
\item[(2.1)] If $B\equiv B_0C^-$, then $AB_0$ and $A^-CB_0^-C^- $ are
             reduced words and $AB_0A^-CB_0^-C^- = ABCA^-B^-C^-$. In
             addition if $AB_0A^-,CB_0^-C^-\in F_m$, then
             $AB,CA^-B^-C^-\in F_m$. 
\item[(2.2)] If $C\equiv B^-C_0$, then $AC_0$ and $A^-B^-C_0^-B$ are
             reduced words and $AC_0A^-B^-C_0^-B = ABCA^-B^-C^-$. In
             addition if $AC_0A^-,B^-C_0^-B\in F_m$, then
             $AB,CA^-B^-C^-\in F_m$. 
\item[(2.3)] If $B\equiv B_0Z^-$ and $C\equiv ZC_0$ for non-empty words
             $B_0$ and $C_0$ and $B_0C_0$ is
             reduced, then $AB_0C_0A^-ZB_0^-C_0^-Z^-$ is reduced
             and $AB_0C_0A^-ZB_0^-C_0^-Z^- 
             = ABCA^-B^-C^-$. In addition if
             $AB_0C_0A^-,ZB_0^-C_0^-Z^-\in F_m$, then
             $AB,CA^-B^-C^-\in F_m$.  
\end{itemize}
\end{lmm}
\begin{proof}
(2.1) The first proposition is obvious. Let $B_0 \equiv XB_1X^-$ for a
 cyclically reduced word $B_1$. Since $(AX)B_1(AX)^-,(CX)B_1^-(CX)^-
\in F_m$, $AX, CX,B_1\in F_m$ by Lemma~\ref{lmm:reducedform}. Now
 $AB = (AX)B_1(CX)^-\in F_m$ and $CA^-B^-C^- 
= (CX)(AX)^-(CB_0^-C^-) \in F_m$. 
We see (2.2) similarly. 

For (2.3) observe the following. Since the both $B_0$ and $C_0$ are
 non-empty, $B_0C_0$ and $B_0^-C_0^-$ are cyclically reduced. Hence,
 using Lemmas~\ref{lmm:reducedform} and \ref{lmm:transform}, we have
 (2.3). 
\end{proof}
The next two lemmas are straightforward and we omit the proofs.
\begin{lmm}\label{lmm:induction3}
Let $A,B,C\in \mathcal{W}(E_n)$ be reduced words (possibly empty) such
 that $ABA^-CB^-C^- \neq e$ and $AB$ and $A^-CB^-C^- $ are reduced. 
Then the following hold: 
\begin{itemize}
\item[(3.1)] If $A\equiv A_0B$, then $A_0B$ and $A_0^-CB^-C^-$ are
             reduced and $A_0BA_0^-CB^-C^- = ABA^-CB^-C^-$. In addition
             if $A_0BA_0^-,CB^-C^-\in F_m$, then $ABA^-,CB^-C^-
             \in F_m$. 
\item[(3.2)] If $B\equiv B_0A$, then $AB_0$ and $CA^-B_0^-C^-$ are
             reduced and $AB_0CA^-B_0^-C^- = ABA^-CB^-C^-$. In addition
             if $AB_0,CA^-B_0^-C^-\in F_m$, then $ABA^-, CB^-C^-\in F_m$. 
\item[(3.3)] If $B\equiv B_0Z$ and $A\equiv A_0Z$ for non-empty words
             $A_0$ and $B_0$ and $B_0A_0^-$ is reduced, then 
             $A_0ZB_0A_0^-CZ^-B_0^-C^-$ is reduced. In addition if
             $A_0ZB_0A_0^-, CZ^-B_0^-C^-\in F_m$, then
             $ABA^-,CB^-C^-\in F_m$. 
\end{itemize}
\end{lmm}
\begin{lmm}\label{lmm:induction4}
Let $A,B,C\in \mathcal{W}(E_n)$ be reduced words (possibly empty) such
 that $ABA^-CB^-C^- \neq e$ and $A$ and $BA^-CB^-C^- $ are reduced
 words. 
Then the following hold: 
\begin{itemize}
\item[(4.1)] If $A\equiv A_0B^-$, $A_0$ and $BA_0^-CB^-C^- $ are
             reduced and $A_0BA_0^-CB^-C^- = ABA^-CB^-C^-$. In addition
             if $A_0BA_0^-,CB^-C^-\in F_m$, then $ABA^-,CB^-C^-\in F_m$. 
\item[(4.2)] If $B\equiv A^-B_0$, and $B_0A^-CB_0^-AC^-$ is reduced and
             $B_0A^-CB_0AC^- = ABA^-CB^-C^- $. In addition if $B_0A^-, 
             CB_0^-AC^-\in F_m$, then $ABA^-, CB^-C^-\in F_m$. 
\item[(4.3)] If $A\equiv A_0Z^-$ and $B\equiv ZB_0$ for non-empty words
             $A_0$, $B_0$ and $A_0B_0$ is reduced, then 
             $A_0B_0ZA_0^-CB_0^-Z^-C^-$ is reduced and
             $A_0B_0ZA_0^-CB_0^-Z^-C^- = ABA^-CB^-C^-$. In addition if 
             $A_0B_0ZA_0^-,$ $CB_0^-Z^-C^-\in F_m$, then
             $ABA^-,CB^-C^-\in F_m$.  
\end{itemize}
\end{lmm}
\begin{lmm}\label{lmm:crucial1}
Let $A,B,C,D\in \mathcal{W}(E_n)$ be reduced non-empty words. 
\begin{itemize}
\item[(1)] if $ABA^-B^-$ is reduced and $ABA^-B^-\in F_m$ and at least
           one of $A, B$ is not small, then $A, B \in F_m$; 
\item[(2)] if $ABCA^-B^-C^-$ is reduced and $ABCA^-B^-C^-\in F_m$ at least
           one of $A,B,C$ is not small, then $A, B, C \in F_m$; 
\item[(3)] if $CABC^-DA^-B^-D^-$ is reduced and $CABC^-DA^-B^-D^-\in F_m$, then
           $A, B, C, D \in F_m$.
\item[(4)] if $CAC^-DA^-D^-$ is reduced and $CAC^-DA^-D^-\in F_m$, then
           $CAC^-, DA^-D^- \in F_m$.
\end{itemize}
\end{lmm}
\begin{proof}
The statements (1) and (2) are paraphrases of Lemma~\ref{lmm:three}. 

\noindent
(3) Let $c$ be the head of $C$ and $d$ be the tail of $D^-$. Since $c^-$
 and $d^-$ are contiguous, we have $CABC^-, DA^-B^-D^-\in F_m$. 
Since $AB$ and $A^-B^-$ are reduced and the both $A$ and $B$ are
 non-empty, $AB$ is cyclically reduced. Now the conclusion follows from
 Lemmas~\ref{lmm:reducedform} and \ref{lmm:transform}. 

\noindent
(4) This follows from a reasoning in the proof of (3). 
\end{proof}
\begin{lmm}\label{lmm:crucial2}
Let $A^-B^-$ and $X_0ABX_0^-$ be reduced words such that $X_0AB\equiv BAX_1$
 for some $X_1$. If\, $lh(X_0)\le lh(B)$, then there exist $A',B'$ 
such that $lh(B')< lh(B)$, $(A')^-(B')^-$ and $X_0A'B'X_0^-$ are reduced 
words, $X_0A'B'\equiv B'A'X_1$, $A^-B^-X_0ABX_0^- = 
(A')^-(B')^-X_0A'B'X_0^-$, and $A,B\in \langle X_0,A',B'\rangle$.
\end{lmm}
\begin{proof}
First we remark that $lh(X_0)\neq lh(B)$ since $BX_0^-$ is reduced. 
Hence $lh(B)>lh(X_0)$. 
If $lh(B) = lh(X_0)+lh(A)$, then we have $X_0A\equiv B \equiv AX_1$ and
 have the conclusion, i,e, $A'\equiv A$ and $B'\equiv \emptyset$. 

If $lh(B) < lh(X_0)+lh(A)$, we have 
$k>0$ and $A_0,A_1$ such that $B\equiv X_0A_0A_1$, $A\equiv (A_0A_1)^kA_0$, 
 and $A_1$ is non-empty. (We remark that $A_0$ may be empty.) 
Let $A'\equiv A_0$ and $B'\equiv A_1$. 
Since $lh(X_0)+lh(A) = lh(B) + (k-1)lh(A_0A_1)+lh(A_0)$, we have
 $B\equiv A_1A_0X_1$. Let $A'\equiv A_0$ and $B'\equiv A_1$, then we have the
 conclusion. 

If $lh(B) > lh(X_0)+lh(A)$, we have 
$k>0$ and $B_0,B_1$ such that $B_0B_1\equiv X_0A$, $B\equiv (B_0B_1)^kB_0$, 
 and $B_1$ is non-empty. (We remark that $B_0$ may be empty.) 
Since $lh(B_1B_0)=lh(AX_1)$, we have $B_1B_0\equiv AX_1$. Now 
$B\equiv X_0A(B_0B_1)^{k-1}B_0 \equiv (B_0B_1)^{k-1}B_0AX_1$ holds. 
Let $A'\equiv A$ and $B'\equiv (B_0B_1)^{k-1}$, then we have the conclusion.
\end{proof}
In Lemma~\ref{lmm:crucial2} we have $A^-B^-X_0ABX_0^- =
X_1X_0^- = (A')^-(B')^-X_0A'B'X_0^-$. 
\begin{lmm}\label{lmm:remaining}
Let $A,B,X,Y\in \mathcal{W}(E_n)$ be reduced words (possibly empty) such
 that $X$ and $Y$ are non-empty, $Y^-A^-B^-YX^-ABX \neq e$, $Y^-A^-B^-Y$
 and $X^-ABX$ are reduced words, and the reduced word of
 $Y^-A^-B^-YX^-ABX$ is cyclically reduced. 

If $Y^-A^-B^-YX^-ABX \in F_m$, then 
\begin{itemize}
\item[(1)] $Y^-A^-B^-Y, X^-ABX\in F_m$, or 
\item[(2)] $Y^-A^-B^-YX^-ABX$ is equal to $c_s$ or $c_s^-$ for some $s$
 such that $lh(s) = m$ and $s$ is binary branched. 
\end{itemize}
\end{lmm}
\begin{proof}
If $YX^-$ is reduced, then $Y^-A^-B^-YX^-ABX$ is cyclically reduced. By
 an argument analyzing the head and the tail of $Y^-$ and $X$ we can see 
$Y^-A^-B^-Y, X^-ABX\in F_m$. 

Otherwise, in the cancellation of $Y^-A^-B^-YX^-ABX$ the leftmost $Y^-$
 or the rightmost $X$ is deleted. Since $Y^-A^-B^-YX^-ABX \neq e$ and
 $lh(Y^-A^-B^-Y) = 2lh(Y)+lh(AB)$ and $lh(X^-ABX) = 2lh(X)+lh(AB)$,
 $lh(X) \neq lh(Y)$. 
We suppose that $lh(X)>lh(Y)$, i.e. the head of $Y^-$ is deleted. Then
 we have $X \equiv ZY$ for a non-empty word $Z$. 

We first analyze a reduced word of $A^-B^-Z^-ABZ$, where $A^-B^-$ is deleted. 
The head part of $Z^-AB$ is $BA$. 
Applying Lemma~\ref{lmm:crucial2} for $X_0\equiv Z^-$ and $X_1$ repeatedly, we
 have reduced words $A_0$ and $B_0$ such that 
$Z^-A_0B_0Z$ is reduced, $Z^-A_0B_0\equiv B_0A_0X_1$ for some $X_1$,
 $A_0^-B_0^-Z^-A_0B_0Z = A^-B^-Z^-ABZ$, 
 $A,B\in\langle Z, A_0, B_0\rangle$ and $lh(B_0)<lh(Z)$. 

It never occurs that the both $A_0$ and $B_0$ are empty, but one of
 $A_0$ and $B_0$ may be empty. If $B_0=\emptyset$, interchange the role
 of $A_0$ and $B_0$ and by Lemma~\ref{lmm:crucial2} we can assume $B_0$
 is non-empty and $lh(B_0)<lh(Z)$.

First we deal with the case $A_0$ is empty. 
Since the left most $B_0^-$ is deleted in the reduction of
 $B_0^-Z^-B_0Z$, we have non-empty $Z_0$ such that $Z \equiv Z_0B_0^-$
 and have a reduced word $Z_0^-B_0Z_0B_0^-$ with $Z_0^-B_0Z_0B_0^- 
= B_0^-Z^-B_0Z$. Since the left most $Y^-$ is deleted in the reduction of
 $Y^-B_0^-Z^-B_0ZY$ and $Z_0^-B_0Z_0B_0^-Y$ is reduced,
 $Z_0^-B_0Z_0B_0^-$ is cyclically reduced and hence the reduced word of
$Y^-A^-B^-YX^-ABX$ is a cyclical transformation of $Z_0^-B_0Z_0B_0^-$. 
By the fact that $Y$ is the head part of $B_0^-Z^-B_0ZY$, $Y$ is of the
 form $(Z_0^-B_0Z_0B_0^-)^kY_0$ where $Y_0Y_1\equiv Z_0^-B_0Z_0B_0^-$
 for some non-empty $Y_1$ and $k\ge 0$. 

If $Y_0$ is empty, we have $Y^-A^-B^-YX^-ABX = Z_0^-B_0Z_0B_0^-$. If one
 of $Z_0$ and $B_0$ is not small, then $Z_0,B_0\in F_m$ by
 Lemma~\ref{lmm:three} and we have $Y^-A^-B^-Y, X^-ABX\in F_m$ by
 Lemma~\ref{lmm:crucial2} and the fact
 $Y=(Z_0^-B_0Z_0B_0^-)^k$. Otherwise, i.e., when 
 of $Z_0$ and $B_0$ are small, $Y^-A^-B^-YX^-ABX = Z_0^-B_0Z_0B_0^-$ 
is equal to $c_s$ or $c_s^-$ for some $s$ such that $lh(s) = m$ and
$s$ is binary branched by Lemma~\ref{lmm:empty}.

If $Y_0\equiv Z_0^-$, $Y_0\equiv Z_0^-B_0$ or 
$Y_0\equiv Z_0^-B_0Z_0$, the argument is similar to the case that $Y_0$
 is empty. Otherwise $Y_0$ cut short $Z_0^-$, $B_0$, $Z_0$ or
 $B_0^-$. Since arguments are similar, we only deal with the case that
 $Y_0\equiv Z_0^-B_1$ where $B_1B_2\equiv B_0$ for non-empty $B_1$ and
 $B_2$. Then $Y^-A^-B^-YX^-ABX = B_2Z_0B_2^-B_1^-Z_0^-B_1$ and hence 
$B_2Z_0B_2^-, B_1^-Z_0^-B_1\in F_m$ by Lemma~\ref{lmm:crucial1} (4). 
Let $Z_1$ be a cyclically reduced word such that $Z_0\equiv U^-Z_1U$.  
Then $Z_1, B_2U^-, UB_1\in F_m$ by Lemma~\ref{lmm:reducedform}. Now 
\begin{eqnarray*}
Y^-Z_0Y &=& B_1^-Z_0(B_1B_2Z_0^-B_2^-B_1^-Z_0)^k
Z_0(Z_0^-B_1B_2Z_0B_2^-B_1^-)^kZ_0^-B_1 \\
&=& 
(B_1^-Z_0B_1B_2Z_0^-B_2^-)^k
B_1^-Z_0B_1(B_2Z_0B_2^-B_1^-Z_0^-B_1)^k \\
Y^-B_0Y &=& B_1^-Z_0(B_1B_2Z_0^-B_2^-B_1^-Z_0)^k
B_1B_2(Z_0^-B_1B_2Z_0B_2^-B_1^-)^kZ_0^-B_1 \\
&=& B_1^-Z_0B_1(B_2Z_0^-B_2^-B_1^-Z_0B_1)^kB_2Z_0^-B_1
(B_2Z_0B_2^-B_1^-Z_0^-B_1)^k. 
\end{eqnarray*}
Hence $Y^-Z_0Y, Y^-B_0Y\in F_m$. 
Since $Z = Z_0B_0^-$ and $A,B\in \langle Z,B_0\rangle$, we have 
$Y^-ABY, X^-A^-B^-X\in F_m$.   

Next we suppose that $A_0$ is non-empty. 
We have $k>0$ and $A_1$ and $A_2$ such that $Z^-\equiv B_0A_1A_2$,
 $A_0\equiv (A_1A_2)^kA_1$, $X_1\equiv A_2A_1B_0$. Since $X^-AB\equiv UX_1$
 for some $U$ and $X^-ABZ$ is reduced, $X_1Z\equiv A_2A_1B_0A_2^-A_1^-B_0^-$ 
is a reduced word. By the assumption a reduced word of 
$Y^-A_2A_1B_0A_2^-A_1^-B_0^-Y$ 
is cyclically reduced and $A_2A_1B_0A_2^-A_1^-B_0^-Y$ is reduced,
 hence $X_1Z\equiv A_2A_1B_0A_2^-A_1^-B_0^-$ is
 cyclically reduced and the reduced word of
 $Y^-A_2A_1B_0A_2^-A_1^-B_0^-Y$ is given by a cyclical transformation of
 $A_2A_1B_0A_2^-A_1^-B_0^-$. Hence $Y\equiv (A_2A_1B_0A_2^-A_1^-B_0^-)^kY_0$ 
where $k\ge 0$ and $A_2A_1B_0A_2^-A_1^-B_0^-\equiv Y_0Y_1$ for some $Y_1$. 

For instance the reduced word of $Y^-A_2A_1B_0A_2^-A_1^-B_0^-Y$ is of
 the form $B_0A_2^-A_1^-B_0^-A_2A_1$ or $B_2A_2^-A_1^-B_2^-B_1^-A_2A_1B_1$
 where $B_0\equiv B_1B_2$. 
By Lemma~\ref{lmm:crucial1} (4) or (3) respectively we conclude 
 $A_1,A_2,B_0\in F_m$ or $A_1,A_2,B_1,B_2\in F_m$ which implies $Y^-ABY, 
X^-A^-B^-X\in F_m$.   
\end{proof}

\begin{lmm}\label{lmm:main}
For every grope group $G^S$ the following hold: 

If $e\neq [u,v]\in F_m$ and at least one of $u$ and $v$ does not belong
 to $F_m$, then $[u,v]$ is conjugate to $c_s$ or $c_s^-$ in $F_m$ for some $s$ 
 such that $lh(s) = m$ and $s$ is binary branched. 
\end{lmm}

\begin{proof} 
We have $n > m$ such that $u,v\in F_n$. 
It suffices to show the lemma in case that the reduced word for $[u,v]$
 is cyclically reduced. For, suppose that we have the conclusion of the
 lemma in the indicated case. Let $[u,v] \in F_m$ and 
$[u,v] = XYX^-$ where $XYX^-$ is a reduced word and $Y$ is cyclically
 reduced. Then we have $[X^-uX,X^-vX] = X^-[u,v]X = Y$. On the other
 hand $X,Y\in F_m$ by Lemma~\ref{lmm:reducedform}. By the assumption 
at least one of $X^-uX$ and $X^-vX$ does not belong to $F_m$. 
Since $[u,v]$ is conjugate to $Y$ in $F_m$, we have the
 conclusion. 

Let $u,v\in F_n$ such that $[u,v]\neq e$ and the reduced word for
 $[u,v]$ is cyclically reduced. There exist a cyclically
reduced non-empty word $V_0\in \mathcal{W}(E_n)$ and a reduced word
$X\in \mathcal{W}(E_n)$ such that $v = X^-V_0X$ and the
 word $X^-V_0X$ is reduced. Let $U_0$ be a reduced word for $uX^-$. 
Since $V_0$ is a cyclically reduced word, at least one of $U_0V_0$
 and $V_0U_0^-$ is reduced. 
When $U_0V_0$ is reduced, there exist $k\ge 0$ and reduced words
 $Y,A,B$ such that $Y^-ABY$ is reduced, $U_0\equiv Y^-AV_0^k$ 
and $V_0\equiv BA$. 
When $V_0U_0^-$ is reduced, there exist $k\ge 0$ and reduced words
 $Y,A,B$ such that $Y^-ABY$ is reduced, $U_0\equiv Y^-A(V_0^-)^k$ 
and $V_0\equiv BA$. 
In the both bases $uvu^{-1} = Y^-ABY$ and $v = X^-BAX$. 
We remark that $AB$ and $BA$ are cyclically reduced. 

We analyze a reduction procedure of $Y^-ABYX^-A^-B^-X$ in the following. 

\noindent
(Case 0): $X$ and $Y$ are empty. 

In this case the both $A$ and $B$ are non-empty and corresponds to
 Lemma~\ref{lmm:induction1}. Using (1.1) and (1.2) alternately and
 (1.3) possibly as the last step we obtain a reduced word of
 $ABA^-B^-$. If the reduced word $XYZX^-Y^-Z^-$ satisfies that one of
 $X,Y,Z$ is not small, by (1) and (2) of Lemma~\ref{lmm:crucial1} and
 applying  Lemma~\ref{lmm:induction1} repeatedly we can see $A,B\in F_m$. 
Otherwise, one of $X,Y,Z$ is empty and $[u,v] = c_s$ or $[u,v] = c_s^-$
 for some binary branched $s$ with $lh(s) = m$ by Lemma~\ref{lmm:empty}. 

\noindent
(Case 1): $Y$ is empty, but $X$ is non-empty. 

\noindent
(Case 2): $X$ is empty, but $Y$ is non-empty. 

In these cases arguments are symmetric, we only deal with (Case 1). 
There is possibility that one of $A$ and $B$ may be empty, though at
 least one of $A$ and $B$ is non-empty. We assume that $A$ is non-empty. 
We trace Lemmas~\ref{lmm:induction2}, \ref{lmm:induction3}, 
\ref{lmm:induction4} to get a reduced word of $ABX^-A^-B^-X$. Then we
 apply one of (2), (3) and (4) of Lemma~\ref{lmm:crucial1} to the
 reduced word and applying Lemma~\ref{lmm:induction1} repeatedly we get 
a reduced word. Then we have $A,B\in F_m$, which implies $u,v\in F_m$, or 
$[u,v] = c_s$ etc. as in (Case 0). 

\noindent
(Case 3): The both $X$ and $Y$ are non-empty. 

Only in this case we use the assumption that the reduced word of
 $Y^-ABYX^-A^-B^-X$ is cyclically reduced. 
By Lemma~\ref{lmm:remaining} we have the conclusion. 
\end{proof}
\begin{lmm}\label{lmm:free1}
Let $F$ be a free group generated by $C$ and $c,d\in C$ be distinct
 elements. If $[c,d] = [u,v]$ for $u,v\in F$, then neither $u$ nor $v$
 belongs to the commutator subgroup of $F$. 
\end{lmm}
\begin{proof}
Since $c,d$ are generators, $[c,d] \notin [F, [F,F]]$ and the conclusion
 follows. 
\end{proof}
\begin{lmm}\label{lmm:free2}
Let $F$ be a free group generated by $B$ and $b_0,b_1\in B$ be
 distinct. 
If $c,d\in \{ b,b^- : b\in B\}$ and $[b_0,b_1] = [x^{-1}cx, y^{-1}dy]$
 for $x,y\in F$, then $c,d\in \{ b_0, b_0^-,b_1,b_1^-\}$ and moreover
$c \in \{ b_0,b_0^-\}$ iff $d \in \{ b_1,b_1^-\}$ and 
$c \in \{ b_1,b_1^-\}$ iff $d \in \{ b_0,b_0^-\}$. 
\end{lmm}
\begin{proof}
Using a canonical projection to $\langle b_0,b_1\rangle$ we easily see
 that $c,d\in \{ b_0, b_0^-,b_1,b_1^-\}$. To see the remaining part it
 suffices to show that if $c=b_0$, and $d = b_0$ or $b_0^-$, 
then $[b_0,b_1] \neq [x^{-1}cx, y^{-1}dy]$ for any $x,y$. 

We show that $b_0b_1b_0^-b_1^-$ is not cyclically equivalent to 
the reduced word for $[x^{-1}cx, y^{-1}dy]$. For this purpose we may
 assume $x=e$. We only deal with $d =b_0$. 
We have a reduced word $Y$ such that $y^{-1}b_0y = Y^-b_0Y$
 and $Y^-b_0Y$ is reduced. (Note that $y=Y$ may not hold.)
The head of $Y$ is not $b_0$ nor $b_0^-$, since $Y^-b_0Y$ is
 reduced. When the tail of $Y$ is $b_0$ or $b_0^-$, we choose $n\ge 0$
 so that $Y\equiv Zb_0^n$ or $Y\equiv Z(b_0^-)^n$ respectively and $n$ is
 maximal. Then $Z$ is non-empty. Now $b_0Z^-b_0Zb_0^-Z^-b_0^-Z$ is a
 cyclically reduced word which is cyclically equivalent to  
$b_0Y^-b_0Yb_0^-Y^-b_0^-Y$. 
Since $b_0Z^-b_0Zb_0^-Z^-b_0^-Z$ is not cyclically equivalent to
 $b_0b_1b_0^-b_1^-$, we have the conclusion. 
\end{proof}

{\it Proof of\/} Theorem~\ref{thm:main}. 
Let $h:G^{S_0}\to G^S$ be a nontrivial homomorphism. 
Then there exists $s_*\in S_0$ such that $h(c_{s_*})$ is nontrivial
(clearly for every finite sequence $s$ starting with $s_*$ also $h(c_s)$ is nontrivial).
We let $c_s = c^{S_0}_s$ and $d_t = c^{S}_t$ and $F_m = F^S_m$. 

We have $n$ such that $h(c_{s_*}) \in F_n$. Since $F_n$ is
free, $\op{Im}(h)$ is not included in $F_n$ and hence there exists $s_0\in
S_0$ starting with $s_*$ and such that $h(c_{s_0})\in F_n$, but $h(c_{s_00})\notin F_n$ or
$h(c_{s_01})\notin F_n$. Then by Lemma~\ref{lmm:main} we have 
$d_{t_0}\in E_n$ such that $h(c_{s_0})$ is conjugate to $d_{t_0}$ or
$d_{t_0}^-$ and $t_0$ is binary branched. 

Moreover, Lemma~\ref{lmm:antinormal} implies that 
neither $h(c_{s_00})$ nor $h(c_{s_01})$ belongs to $F_n$. 
We show the following by induction on $k\in \mathbb{N}$: 

\smallskip
\noindent
(1) For $u\in Seq(\un{2})$ with $lh(u)=k$ 
\begin{itemize}
\item[(a)]  $h(c_{s_0u})$ is conjugate to $d_{t_0v}$ or $d_{t_0v}^-$ in
            $F_{n+k}$ and $t_0v$ is binary branched for some
            $v\in Seq(\un{2})$ with $lh(v)=k$;  
\item[(b)]  Neither $h(c_{s_0u0})$ nor $h(c_{s_0u1})$ belongs to
            $F_{n+k}$; 
\end{itemize}

\smallskip
\noindent
(2) For every $v\in Seq(\un{2})$ with $lh(v)=k$ there exists $u\in
Seq(\un{2})$ such that $lh(u) = k$ and $h(c_{s_0u})$ is conjugate to
$d_{t_0v}$ or $d_{t_0v}^-$ in $F_{n+k}$. 

We have shown that this holds when $k = 0$. 

Suppose that (1) and (2) hold for $k$. Let $lh(u) = k$ and 
$h(c_{s_0u})$ is conjugate to $d_{t_0v}$ or $d_{t_0v}^-$ etc. 
Then $[h(c_{s_0u0}), h(c_{s_0u1})]$ is conjugate to
$[d_{t_0v0},d_{t_0v1}]$ or $[d_{t_0v1},d_{t_0v0}]$ in 
$F_{n+k+1}$. 
We claim $h(c_{s_0u0})\in F_{n+k+1}$. To show this by
contradiction, suppose that $h(c_{s_0u0})\notin F_{n+k+1}$. 
Apply Lemma~\ref{lmm:main} to $F_{n+k+1}$, then we have
$[h(c_{s_0u0}),h(c_{s_0u1})]$ is a conjugate to $d_t$ or $d_t^-$ with
$lh(t) = n + k + 1$ in $F_{n+k+1}$, which is impossible since
$[h(c_{s_0u0}),h(c_{s_0u1})]\in [F_{n+k+1},F_{n+k+1}]$. 
Similarly we have $h(c_{s_0u1})\in F_{n+k+1}$. 

On the other hand, neither $h(c_{s_0u0})$ nor $h(c_{s_0u1})$
belongs to $[F_{n+k+1},F_{n+k+1}]$ by
Lemma~\ref{lmm:free1}. Hence 
at least one of $h(c_{s_0u00})$ and $h(c_{s_0u01})$
does not belong to $F_{n+k+1}$ and consequently neither $h(c_{s_0u00})$
nor $h(c_{s_0u01})$ belongs to $F_{n+k+1}$ by
Lemma~\ref{lmm:antinormal}.  

Hence $h(c_{s_0u0})$ is conjugate to $d_t$ or $d_t^-$ with $lh(t) =
n + k + 1$ by Lemma~\ref{lmm:main}. Similarly, $h(c_{s_0u1})$ is conjugate
to $d_{t'}$ or $d_{t'}^-$ with $lh(t') = n + k + 1$. 
Since $[h(c_{s_0u0}), h(c_{s_0u1})]$ is conjugate to
$[d_{t_0v0},d_{t_0v1}]$ or $[d_{t_0v0},d_{t_0v1}]$ in $F_{n+k+1}$, 
$h(c_{s_0u0})$ and $h(c_{s_0u1})$ are conjugate to 
$d_{t_0vj}$ or $d_{t_0vj}^-$ for some $j\in \un{2}$ and for each 
$j\in \un{2}$ the element $d_{t0vj}$ is conjugate to exactly one of
$h(c_{s_0u0})$, $h(c_{s_0u1})$, $h(c_{s_0u0})^-$ and $h(c_{s_0u1})^-$ by
Lemma~\ref{lmm:free2}.  
Hence (1) and (2) hold for $k+1$.
Now we have shown the induction step and finished the proof. 
\qed
\begin{rmk}
Though the conclusion of Theorem~\ref{thm:main} is rather simple,
 embeddings from $G^{S_0}$ into $G^S$ may be complicated. In particular
 automorphisms on $G^{S_0}$ may be complicated, since the following
 hold: 
\[
 [dc^-d^-, dcd^-c^-d^-] = dc^-d^-dcd^-c^-d^-dcd^-dcdc^-d^- = cdc^-d^- =
 [c,d]. 
\]

\end{rmk}


\end{document}